ORIGINAL RESEARCH

# An order theoretic approach in fixed point theory

Yaé Ulrich Gaba

**Abstract** In the present article, we show the existence of a coupled fixed point for an order preserving mapping in a preordered left $K$-complete quasi-pseudometric space using a preorder induced by an appropriate function. We also define the concept of left-weakly related mappings on a preordered space and discuss common coupled fixed points for two and three left-weakly related mappings in the same space. Similar results are given for right-weakly related mappings, the dual notion of left-weakly related mappings.

## Introduction

Fixed point theory plays a major role in many applications including variational and linear inequalities, optimization and applications in the field of approximation theory and minimum norm problem. The first important result on fixed points for contractive-type mapping was the well-known Banach's contraction principle that appeared in explicit form in his thesis in 1922, where it was used to establish the existence of a solution for an integral equation. This theorem is a key result in non-linear analysis. Another interesting result on fixed points for contractive-type mapping is due to Edelstein (1962) who actually obtained slightly more general versions. Many years later, another direction of such generalizations (see [1, 2]) has been obtained by weakening the requirements in the contractive condition and in compensation, by simultaneously enriching the metric space structure with a partial order.

In this process of generalization, the study of common fixed points of mappings satisfying certain contractive conditions has also been at the center of rigorous research activity, see [3].

Bhashkar and Lakshmikantham [4] introduced the concept of a coupled fixed point of a mapping $F : X \times X \to X$ (where $X$ a non-empty set) and established some coupled fixed point theorems in partially ordered complete metric spaces. By doing so, they opened the way to a flourishing sub-area in the fixed point theory.

Moreover, in the last few years, there has been a growing interest in the theory of quasi-metric spaces and other related structures such as asymmetric normed linear spaces (see for instance [5]). This theory provides a convenient framework in the study of several problems in theoretical computer science and approximation theory. It is in this setting that we give our results.

The aim of this paper is to analyze the existence of common and coupled fixed points for mapping defined on a left $K$-complete quasi-pseudometric space $(X, d)$. The technique of proof is different and more natural in the sense we do not use any contractive conditions. In our work, we show the existence of a coupled fixed point for an isotone mapping in a preordered left $K$-complete quasi-pseudometric space using a preorder induced by an appropriate function $\phi$. Furthermore, common coupled fixed point for two and three mappings satisfying a certain relation that we specify later is also discussed in the same space.

Y. U. Gaba (✉)
Department of Mathematics and Applied Mathematics,
University of Cape Town, Rondebosch 7701, South Africa
e-mail: gabayae2@gmail.com

## Preliminaries

In this section, we recall some elementary definitions from the asymmetric topology and the order theory, which are necessary for a good understanding of the work below.

**Definition 2.1** Consider a non-empty set $X$ and a binary relation $\preceq$ on $X$. Then, $\preceq$ is a preorder, or quasiorder, if it is reflexive and transitive, i.e., for all $a, b$ and $c \in X$, we have that:

- $a \preceq a$ (reflexivity);
- if $a \preceq b$ and $b \preceq c$, then $a \preceq c$ (transitivity).

A set that is equipped with a preorder is called a preordered space (or proset).

**Definition 2.2** Let $(X, \preceq_X)$ and $(Y, \preceq_Y)$ be two prosets. A map $T : X \to Y$ is said to be preorder-preserving or isotone if for any $x, y \in X$,

$$x \preceq_X y \implies Tx \preceq_Y Ty.$$

Similarly, for any family $(X_i, \preceq_{X_i}), i = 1, 2, \ldots, n;\ (Y, \preceq_Y)$ of posets, a mapping $F : X_1 \times X_2 \times \cdots \times X_n \to Y$ is said to be preorder-preserving or isotone if for any $(x_1, x_2, \ldots, x_n)$, $(z_1, z_2, \ldots, z_n) \in X_1 \times X_2 \times \cdots \times X_n$,

$$x_i \preceq_{X_i} z_i \text{ for all } i = 1, 2, \ldots, n \implies F(x_1, x_2, \ldots, x_n) \preceq_Y F(z_1, z_2, \ldots, z_n).$$

Dually, we have

**Definition 2.3** (*Compare* [6]) Let $X$ be a non-empty set. A function $d : X \times X \to [0, \infty)$ is called a quasi-pseudometric on $X$ if:

(i) $d(x, x) = 0 \quad \forall\, x \in X$,
(ii) $d(x, z) \leq d(x, y) + d(y, z) \quad \forall\, x, y, z \in X$.

Moreover, if $d(x, y) = 0 = d(y, x) \implies x = y$, then $d$ is said to be a $T_0$-quasi-pseudometric. The latter condition is referred to as the $T_0$-condition.

*Remark 2.4*

- Let $d$ be a quasi-pseudometric on $X$, then the map $d^{-1}$ defined by $d^{-1}(x, y) = d(y, x)$ whenever $x, y \in X$ is also a quasi-pseudometric on $X$, called the conjugate of $d$. In the literature, $d^{-1}$ is also denoted $d^t$ or $\bar{d}$.
- It is easy to verify that the function $d^s$ defined by $d^s := d \vee d^{-1}$, i.e., $d^s(x, y) = \max\{d(x, y), d(y, x)\}$ defines a metric on $X$ whenever $d$ is a $T_0$-quasi-pseudometric on $X$.

Let $(X, d)$ be a quasi-pseudometric space. For $x \in X$ and $\varepsilon > 0$,

$$B_d(x, \varepsilon) = \{y \in X : d(x, y) < \varepsilon\}$$

denotes the open $\varepsilon$-ball at $x$. The collection of all such balls yields a base for the topology $\tau(d)$ induced by $d$ on $X$. Hence, for any $A \subseteq X$, we shall, respectively, denote by $\text{int}_{\tau(d)}A$ and $\text{cl}_{\tau(d)}A$ the interior and the closure of the set $A$ with respect to the topology $\tau(d)$.

Similarly, for $x \in X$ and $\varepsilon \geq 0$,

$$C_d(x, \varepsilon) = \{y \in X : d(x, y) \leq \varepsilon\}$$

denotes the closed $\varepsilon$-ball at $x$.

**Definition 2.5** Let $(X, d)$ be a quasi-pseudometric space. The convergence of a sequence $(x_n)$ to $x$ with respect to $\tau(d)$, called $d$-convergence or left-convergence and denoted by $x_n \xrightarrow{d} x$, is defined in the following way

$$x_n \xrightarrow{d} x \iff d(x, x_n) \longrightarrow 0. \quad (1)$$

Similarly, the convergence of a sequence $(x_n)$ to $x$ with respect to $\tau(d^{-1})$, called $d^{-1}$-convergence or right-convergence and denoted by $x_n \xrightarrow{d^{-1}} x$, is defined in the following way

$$x_n \xrightarrow{d^{-1}} x \iff d(x_n, x) \longrightarrow 0. \quad (2)$$

Finally, in a quasi-pseudometric space $(X, d)$, we shall say that a sequence $(x_n)$ $d^s$-converges to $x$ if it is both left and right convergent to $x$, and we denote it as $x_n \xrightarrow{d^s} x$ or $x_n \to x$ when there is no confusion. Hence,

$$x_n \xrightarrow{d^s} x \iff x_n \xrightarrow{d} x \text{ and } x_n \xrightarrow{d^{-1}} x.$$

**Definition 2.6** A sequence $(x_n)$ in a quasi-pseudometric $(X, d)$ is called

(a) left $d$-Cauchy if for every $\epsilon > 0$, there exist $x \in X$ and $n_0 \in \mathbb{N}$ such that

$$\forall\, n \geq n_0 \quad d(x, x_n) < \epsilon;$$

(b) left $K$-Cauchy if for every $\epsilon > 0$, there exists $n_0 \in \mathbb{N}$ such that

$$\forall\, n, k : n_0 \leq k \leq n \quad d(x_k, x_n) < \epsilon;$$

(c) $d^s$-Cauchy if for every $\epsilon > 0$, there exists $n_0 \in \mathbb{N}$ such that

$$\forall n, k \geq n_0 \quad d^s(x_n, x_k) < \epsilon.$$

Dually, we define in the same way, right $d$-Cauchy and right $K$-Cauchy sequences.

*Remark 2.7*

- $d^s$-Cauchy $\implies$ left $K$-Cauchy $\implies$ left $d$-Cauchy. The same implications hold for the corresponding right notions. None of the above implications is reversible.

- A sequence is left $d$-Cauchy with respect to $d$ if and only if it is right $K$-Cauchy with respect to $d^{-1}$.
- A sequence is left $K$-Cauchy with respect to $d$ if and only if it is right $K$-Cauchy with respect to $d^{-1}$.
- A sequence is $d^s$-Cauchy if and only if it is both left and right $d$-Cauchy.

**Definition 2.8** (*Compare* [6]) A quasi-pseudometric space $(X, d)$ is called

- left-$K$-complete provided that any left $K$-Cauchy sequence is $d$-convergent,
- left Smyth sequentially complete if any left $K$-Cauchy sequence is $d^s$-convergent.

The dual notions of right-completeness are easily derived from the above definition.

**Definition 2.9** A $T_0$-quasi-pseudometric space $(X, d)$ is called bicomplete provided that the metric $d^s$ on $X$ is complete.

**Definition 2.10** Let $(X, d)$ be a quasi-pseudometric type space. A function $T : X \to X$ is called $d$-sequentially continuous or left-sequentially continuous if for any $d$-convergent sequence $(x_n)$ with $x_n \xrightarrow{d} x$, the sequence $(Tx_n)$ $d$-converges to $Tx$, i.e., $Tx_n \xrightarrow{d} Tx$.

Similarly, a function $T : X \times X \to X$ is said to be $d$-sequentially continuous or left-sequentially continuous if for any sequences $(x_n)$ and $(y_n)$ such that $x_n \xrightarrow{d} x$ and $y_n \xrightarrow{d} y$, then $F(x_n, y_n) \xrightarrow{d} F(x, y)$.

Similarly, we define a $d^{-1}$-sequentially continuous or right-sequentially continuous and $d^s$-sequentially continuous functions.

**Definition 2.11** (*Compare* [2]) Let $X$ be a non-empty set. An element $(x, y) \in X \times X$ is called:

(E1) a coupled fixed point of the mappings $F : X \times X \to X$ if $F(x, y) = x$ and $F(y, x) = y$;
(E2) a coupled coincidence point of the mappings $F : X \times X \to X$ and $T : X \to X$ if $T(x, y) = Tx$ and $T(y, x) = Ty$, in this case $(Tx, Ty)$ is called the coupled point of coincidence;
(E3) a common coupled fixed point of the mappings $F : X \times X \to X$ and $T : X \to X$ if $F(x, y) = Tx = x$ and $F(y, x) = Ty = y$.

**Definition 2.12** Let $X$ be a non-empty set. An element $(x, y) \in X \times X$ is called:

(D1) a common coupled coincidence point of the mappings $F : X \times X \to X$ and $T, R : X \to X$ if $F(x, y) = Tx = Rx$ and $F(y, x) = Ty = Ry$;
(D2) a common coupled fixed point of the mappings $F : X \times X \to X$ and $T, R : X \to X$ if $F(x, y) = Tx = Rx = x$ and $F(y, x) = Ty = Ry = y$.

## First results

We start by the following lemma.

**Lemma 3.1** Let $(X, d)$ be a quasi-pseudometric space and $\phi : X \to \mathbb{R}$ a map. Define the binary relation $"\preceq"$ on $X$ as follows:

$$x \preceq y \iff d(x, y) \leq \phi(y) - \phi(x).$$

Then, $"\preceq"$ is a preorder on $X$. It will be called the preorder induced by $\phi$.

*Proof 3.2*

1. • Reflexivity: For all $x \in X$, $d(x, x) = 0 = \phi(x) - \phi(x)$ hence $x \preceq x$, i.e., $"\preceq"$ is reflexive.
2. • Transitivity: For $x, y, z \in X$ s.t. $x \preceq y$ and $y \preceq z$, we have

$$d(x, y) \leq \phi(y) - \phi(x),$$

and

$$d(y, z) \leq \phi(z) - \phi(y),$$

and since

$$\begin{aligned} d(x, z) &\leq d(x, y) + d(y, z) \\ &\leq \phi(y) - \phi(x) + \phi(z) - \phi(y) \\ &= \phi(z) - \phi(x). \end{aligned}$$

we have $x \preceq z$. Thus, $"\preceq"$ is transitive, and so the relation $"\preceq"$ is a preorder on $X$.

*Remark 3.3* If in addition, the space $(X, d)$ is $T_0$, then the relation $\overline{\preceq}$ defined by

$$x \overline{\preceq} y \iff d^s(x, y) \leq \phi(y) - \phi(x).$$

is a partial order on $X$.

Now, we prove the following theorem.

**Theorem 3.4** *Let $(X, d)$ be a Hausdorff left K-complete $T_0$-quasi-pseudometric space, $\phi : X \to \mathbb{R}$ be a bounded from above function and $"\preceq"$ the preorder induced by $\phi$. Let $F : X \times X \to X$ be a preorder-preserving and $d$-sequentially continuous mapping on $X$ such that there exist two elements $x_0, y_0 \in X$ with*

$$x_0 \preceq F(x_0, y_0) \quad \text{and} \quad y_0 \preceq F(y_0, x_0).$$

*Then, $F$ has a coupled fixed point in $X$.*

*Proof 3.5* Let $x_0, y_0 \in X$ with $x_0 \preceq F(x_0, y_0)$ and $y_0 \preceq F(y_0, x_0)$. We construct the sequences $(x_n)$ and $(y_n)$ in $X$ as follows:

$$x_{n+1} = F(x_n, y_n) \text{ and } y_{n+1} = F(y_n, x_n) \quad \text{for all } n \geq 0. \quad (3)$$

We shall show that

$$x_n \preceq x_{n+1} \quad \text{for all } n \geq 0, \quad (4)$$

and

$$y_n \preceq y_{n+1} \quad \text{for all } n \geq 0. \quad (5)$$

For this purpose, we use the mathematical induction.

Since $x_0 \preceq F(x_0, y_0)$ and $y_0 \preceq F(y_0, x_0)$ and as $x_1 = F(x_0, y_0)$ and $y_1 = F(y_0, x_0)$, we have $x_0 \preceq x_1$ and $y_0 \preceq y_1$. Thus, (4) and (5) hold for $n = 0$.

Suppose that (4) and (10) hold for some $k > 0$. Then, since $x_k \preceq x_{k+1}$ and $y_k \preceq y_{k+1}$ and $F$ is preorder preserving, we have

$$x_{k+1} = F(x_k, y_k) \preceq F(x_{k+1}, y_{k+1}) = x_{k+2} \quad (6)$$

and

$$y_{k+1} = F(y_k, x_k) \preceq F(y_{k+1}, x_{k+1}) = y_{k+2}. \quad (7)$$

Thus, by mathematical induction, we conclude that (4) and (5) hold for all $n \geq 0$. Therefore,

$$x_0 \preceq x_1 \preceq x_2 \preceq \cdots \preceq x_n \preceq \cdots,$$

and

$$y_0 \preceq y_1 \preceq y_2 \preceq \cdots \preceq y_n \preceq \cdots.$$

By definition of the preorder, we have

$$\phi(x_0) \leq \phi(x_1) \leq \phi(x_2) \leq \cdots \leq \phi(x_n) \leq \cdots,$$

and

$$\phi(y_0) \leq \phi(y_1) \leq \phi(y_2) \leq \cdots \leq \phi(y_n) \leq \cdots.$$

Hence, the sequences $(\phi(x_n))$ and $(\phi(y_n))$ are non-decreasing sequences of real numbers. Since $\phi$ is bounded from above, the sequences $(\phi(x_n))$ and $(\phi(y_n))$ are convergent and, therefore, Cauchy. This entails that for any $\varepsilon > 0$, there exists $n_0 \in \mathbb{N}$ such that for any $m > n > n_0$, we have $\phi(x_m) - \phi(x_n) < \varepsilon$ and $\phi(y_m) - \phi(y_n) < \varepsilon$. Since whenever $m > n > n_0$, $x_n \preceq x_m$ and $y_n \preceq y_m$, it follows that

$$d(x_n, x_m) \leq \phi(x_m) - \phi(x_n) < \varepsilon,$$

and

$$d(y_n, y_m) \leq \phi(y_m) - \phi(y_n) < \varepsilon.$$

We conclude that $(x_n)$ and $(y_n)$ are left $K$-Cauchy in $X$ and since $X$ is left $K$-complete, there exist $x^*, y^* \in X$ such that $x_n \xrightarrow{d} x^*$ and $y_n \xrightarrow{d} y^*$. Since $F$ is $d$-sequentially continuous, we have

$$x_n \xrightarrow{d} x^* \iff x_n = F(x_{n-1}, y_{n-1}) \xrightarrow{d} x^* \iff F(x^*, y^*) = x^*,$$

and

$$y_n \xrightarrow{d} y^* \iff y_n = F(y_{n-1}, x_{n-1}) \xrightarrow{d} y^* \iff F(y^*, x^*) = y^*.$$

Thus, we have proved that $F(x^*, y^*) = x^*$ and $F(y^*, x^*) = y^*$, i.e., $(x^*, y^*)$ is a coupled fixed point of $F$.

**Corollary 3.6** Let $(X, d)$ be a Hausdorff right $K$-complete $T_0$-quasi-pseudometric space, $\phi : X \to \mathbb{R}$ be a bounded from below function and $"\preceq"$ the preorder induced by $\phi$. Let $F : X \times X \to X$ be a preorder-preserving and $d^{-1}$-sequentially continuous mapping on $X$ such that there exist two elements $x_0, y_0 \in X$ with

$$x_0 \succeq F(x_0, y_0) \quad \text{and } y_0 \succeq F(y_0, x_0).$$

Then, $F$ has a coupled fixed point in $X$.

**Corollary 3.7** Let $(X, d)$ be a bicomplete $T_0$-quasi-pseudometric space, $\phi : X \to \mathbb{R}$ be a bounded from below function and $"\preceq"$ the preorder induced by $\phi$. Let $F : X \times X \to X$ be a preorder-preserving and $d^s$-sequentially continuous mapping on $X$ such that there exist two elements $x_0, y_0 \in X$ with

$$x_0 \preceq F(x_0, y_0) \quad \text{and } y_0 \preceq F(y_0, x_0).$$

Then, $F$ has a coupled fixed point in $X$.

## Common coupled fixed point

Now, we define the concept of weakly related mappings on preordered spaces as follows:

**Definition 4.1** Let $(X, \preceq)$ be a preordered space, and $F : X \times X \to X$ and $g : X \to X$ be two mappings. Then, the pair $\{F, g\}$ is said to be weakly left-related if the two following conditions are satisfied:

(C1) $F(x, y) \preceq gF(x, y)$ and $gx \preceq F(gx, gy)$ for all $(x, y) \in X \times X$,

(C2) $F(y, x) \preceq gF(y, x)$ and $gy \preceq F(gy, gx)$ for all $(x, y) \in X \times X$.

**Definition 4.2** Let $(X, \preceq)$ be a preordered space, and $F : X \times X \to X$ and $g : X \to X$ be two mappings. Then, the pair $\{F, g\}$ is said to be weakly right-related if the two following conditions are satisfied:

(D1) $gF(x, y) \preceq F(x, y)$ and $F(gx, gy) \preceq gx$ for all $(x, y) \in X \times X$,

(D2) $gF(y,x) \preceq F(y,x)$ and $F(gy, gx) \preceq gy$ for all $(x, y) \in X \times X$.

We now state and prove the first common coupled fixed point existence theorem for the weakly related mappings.

**Theorem 4.3** *Let $(X, d)$ be a Hausdorff left K-complete $T_0$-quasi-pseudometric space, $\phi : X \to \mathbb{R}$ be a bounded from above function and $"\preceq"$ the preorder induced by $\phi$. Let $F : X \times X \to X$ and $G : X \to X$ be two d-sequentially continuous mapping on X such that the pair $\{F, G\}$ is weakly left-related. If there exist two elements $x_0, y_0 \in X$ with*

$$x_0 \preceq F(x_0, y_0) \quad \text{and } y_0 \preceq F(y_0, x_0).$$

*Then, F and G have a common coupled fixed point in X.*

*Proof 4.4* Let $x_0, y_0 \in X$ with $x_0 \preceq F(x_0, y_0)$ and $y_0 \preceq F(y_0, x_0)$. We construct the sequences $(x_n)$ and $(y_n)$ in X as follows:

$$x_{2n+1} = F(x_{2n}, y_{2n}) \text{ and } x_{2n+2} = Gx_{2n+1} \text{ for all } n \geq 0, \quad (8)$$

and

$$y_{2n+1} = F(y_{2n}, x_{2n}) \text{ and } y_{2n+2} = Gy_{2n+1} \text{ for all } n \geq 0. \quad (9)$$

We shall show that

$$x_n \preceq x_{n+1} \quad \text{for all } n \geq 0, \quad (10)$$

and

$$y_n \preceq y_{n+1} \quad \text{for all } n \geq 0. \quad (11)$$

Since $x_0 \preceq F(x_0, y_0)$, using (8), we have and $x_0 \preceq x_1$. Again since the pair $\{F, G\}$ is weakly left-related, we have, from (8) $x_1 = F(x_0, y_0) \preceq GF(x_0, y_0) = Gx_1 = x_2$, i.e., $x_1 \preceq x_2$. Also, since $Gx_1 \preceq F(Gx_1, Gy_1)$, and using (8), we have $x_2 = Gx_1 \preceq F(Gx_1, Gy_1) = F(x_2, y_2) = x_3$, i.e., $x_2 \preceq x_3$. Similarly, using the fact that the pair $\{F, G\}$ is weakly left-related and the relations (8), we get

$$x_0 \preceq x_1 \preceq x_2 \preceq \cdots \preceq x_n \preceq \cdots.$$

A similar reasoning, using fact that the pair $\{F, G\}$ is weakly left-related and the relations (9), leads to

$$y_0 \preceq y_1 \preceq y_2 \preceq \cdots \preceq y_n \preceq \cdots.$$

By definition of the preorder, we have

$$\phi(x_0) \leq \phi(x_1) \leq \phi(x_2) \leq \cdots \leq \phi(x_n) \leq \cdots$$

and

$$\phi(y_0) \leq \phi(y_1) \leq \phi(y_2) \leq \cdots \leq \phi(y_n) \leq \cdots.$$

Hence, the sequences $(\phi(x_n))$ and $(\phi(y_n))$ are non-decreasing sequences of real numbers. Since $\phi$ is bounded from above, the sequences $(\phi(x_n))$ and $(\phi(y_n))$ are convergent and, therefore, Cauchy. This entails that for any $\varepsilon > 0$, there exists $n_0 \in \mathbb{N}$ such that for any $m > n > n_0$, we have $\phi(x_m) - \phi(x_n) < \varepsilon$ and $\phi(y_m) - \phi(y_n) < \varepsilon$. Since whenever $m > n > n_0$, $x_n \preceq x_m$ and $y_n \preceq y_m$, it follows that

$$d(x_n, x_m) \leq \phi(x_m) - \phi(x_n) < \varepsilon,$$

and

$$d(y_n, y_m) \leq \phi(y_m) - \phi(y_n) < \varepsilon.$$

It follows that $(x_n)$ and $(y_n)$ are left K-Cauchy in X and since X is left K-complete, there exist $x^*, y^* \in X$ such that $x_n \xrightarrow{d} x^*$ and $y_n \xrightarrow{d} y^*$.

Since F and G are d-sequentially continuous, it is easy to see, using (8), that

$$x_{2n+1} \xrightarrow{d} x^* \iff x_{2n+1} = F(x_{2n}, y_{2n}) \xrightarrow{d} x^* \iff F(x^*, y^*) = x^*,$$

and

$$x_{2n+2} \xrightarrow{d} x^* \iff x_{2n+2} = Gx_{2n+1} \xrightarrow{d} x^* \iff Gx^* = x^*,$$

and hence

$$Gx^* = x^* = F(x^*, y^*).$$

Similarly, since F and G are d-sequentially continuous, using (9), we easily derive that

$$Gy^* = y^* = F(y^*, x^*).$$

Hence, $(x^*, y^*)$ is a coupled common fixed point of F and G.

**Corollary 4.5** *Let $(X, d)$ be a Hausdorff right K-complete $T_0$-quasi-pseudometric space, $\phi : X \to \mathbb{R}$ be a bounded from below function and $"\preceq"$ the preorder induced by $\phi$. Let $F : X \times X \to X$ and $G : X \to X$ be two $d^{-1}$-sequentially continuous mapping on X such that the pair $\{F, G\}$ is weakly right-related. If there exist two elements $x_0, y_0 \in X$ with*

$$F(x_0, y_0) \preceq x_0 \quad \text{and } F(y_0, x_0) \preceq y_0.$$

*Then, F and G have a common coupled fixed point in X.*

**Corollary 4.6** *Let $(X, d)$ be a bicomplete $T_0$-quasi-pseudometric space, $\phi : X \to \mathbb{R}$ be a bounded from above function and $"\overline{\preceq}"$ the preorder induced by $\phi$. Let $F : X \times X \to X$ and $G : X \to X$ be two $d^s$-sequentially continuous mapping on X such that the pair $\{F, G\}$ satisfies the following two conditions:*

(C1') $F(x, y) \overline{\preceq} GF(x, y)$ and $Gx \overline{\preceq} F(Gx, Gy)$ for all $(x, y) \in X \times X$,
(C2') $F(y, x) \overline{\preceq} GF(y, x)$ and $Gy \overline{\preceq} F(Gy, Gx)$ for all $(x, y) \in X \times X$.

Springer

If there exist two elements $x_0, y_0 \in X$ with

$$x_0 \preceq F(x_0, y_0) \quad \text{and} \quad y_0 \preceq F(y_0, x_0).$$

Then, $F$ and $G$ have a common coupled fixed point in $X$.

**Theorem 4.7** *Let $(X, d)$ be a Hausdorff left $K$-complete $T_0$-quasi-pseudometric space, $\phi : X \to \mathbb{R}$ be a bounded from above function and $"\preceq"$ the preorder induced by $\phi$. Let $F : X \times X \to X$ and $G, H : X \to X$ be three $d$-sequentially continuous mapping on $X$ such that the pairs $\{F, G\}$ and $\{F, H\}$ are weakly left-related. Then, $F, G$ and $H$ have a common coupled fixed point in $X$.*

*Proof 4.8* Let $x_0, y_0 \in X$. We construct the sequences $(x_n)$ and $(y_n)$ in $X$ as follows:

$$Hx_{3n-3} = x_{3n-2}, \; x_{3n-1} = F(x_{3n-2}, y_{3n-2}) \text{ and } x_{3n} = Gx_{3n-1} \quad \text{for all } n \geq 1, \tag{12}$$

and

$$Hy_{3n-3} = y_{3n-2}, \; y_{3n-1} = F(y_{3n-2}, x_{3n-2}) \text{ and } y_{3n} = Gy_{3n-1} \quad \text{for all } n \geq 1. \tag{13}$$

We shall show that

$$x_n \preceq x_{n+1} \quad \text{for all } n \geq 0, \tag{14}$$

and

$$y_n \preceq y_{n+1} \quad \text{for all } n \geq 0. \tag{15}$$

We have $x_1 = Hx_0$. Since the pair $\{F, H\}$ is weakly left-related, we have $x_1 = Hx_0 \preceq F(Hx_0, Hy_0) = F(x_1, y_1) = x_2$. Again since the pair $\{F, G\}$ is weakly left-related, we have $x_2 = F(x_1, y_1) \preceq GF(x_1, y_1) = Gx_2 = x_3$. Similarly, using the fact that the pairs $\{F, G\}$ and $\{F, H\}$ are weakly left-related and the relations (12), we get

$$x_0 \preceq x_1 \preceq x_2 \preceq x_3 \preceq \cdots \preceq x_n \preceq \cdots.$$

A similar reasoning, using fact that the pair $\{F, g\}$ is weakly left-related and the relations (13), leads to

$$y_0 \preceq y_1 \preceq y_2 \preceq y_3 \preceq \cdots \preceq y_n \preceq \cdots.$$

By definition of the preorder, we have

$$\phi(x_0) \leq \phi(x_1) \leq \phi(x_2) \leq \cdots \leq \phi(x_n) \leq \cdots$$

and

$$\phi(y_0) \leq \phi(y_1) \leq \phi(y_2) \leq \cdots \leq \phi(y_n) \leq \cdots.$$

Hence, the sequences $(\phi(x_n))$ and $(\phi(y_n))$ are non-decreasing sequences of real numbers. Since $\phi$ is bounded from above, the sequences $(\phi(x_n))$ and $(\phi(y_n))$ are convergent and, therefore, Cauchy. This entails that for any $\varepsilon > 0$, there exists $n_0 \in \mathbb{N}$ such that for any $m > n > n_0$, we have $\phi(x_m) - \phi(x_n) < \varepsilon$ and $\phi(y_m) - \phi(y_n) < \varepsilon$. Since whenever $m > n > n_0$, $x_n \preceq x_m$ and $y_n \preceq y_m$, it follows that

$$d(x_n, x_m) \leq \phi(x_m) - \phi(x_n) < \varepsilon,$$

and

$$d(y_n, y_m) \leq \phi(y_m) - \phi(y_n) < \varepsilon.$$

It follows that $(x_n)$ and $(y_n)$ are left $K$-Cauchy in $X$ and since $X$ is left $K$-complete, there exist $x^*, y^* \in X$ such that $x_n \xrightarrow{d} x^*$ and $y_n \xrightarrow{d} y^*$.

Since $F, G$ and $H$ are $d$-sequentially continuous, it is easy to see, using (12), that

$$x_{3n-1} \xrightarrow{d} x^* \iff x_{3n-1} = F(x_{3n-2}, y_{3n-2}) \xrightarrow{d} x^*$$
$$\iff F(x^*, y^*) = x^*,$$

and

$$x_{3n} \xrightarrow{d} x^* \iff x_{3n} = Gx_{3n-1} \xrightarrow{d} x^* \iff Gx^* = x^*,$$

also

$$x_{3n-2} \xrightarrow{d} x^* \iff x_{3n-2} = Hx_{3n-3} \xrightarrow{d} x^* \iff Hx^* = x^*,$$

and hence

$$Hx^* = Gx^* = x^* = F(x^*, y^*).$$

Similarly, since $F, G$ and $H$ are $d$-sequentially continuous, using (13), we easily derive that

$$Hy^* = Gy^* = y^* = F(y^*, x^*).$$

Hence, $(x^*, y^*)$ is a coupled common fixed point of $F, G$ and $H$.

**Corollary 4.9** *Let $(X, d)$ be a Hausdorff right $K$-complete $T_0$-quasi-pseudometric space, $\phi : X \to \mathbb{R}$ be a bounded from below function and $"\preceq"$ the preorder induced by $\phi$. Let $F : X \times X \to X$ and $G, H : X \to X$ be three $d^{-1}$-sequentially continuous mapping on $X$ such that the pairs $\{F, G\}$ and $\{F, H\}$ are weakly right-related. Then, $F, G$ and $H$ have a common coupled fixed point in $X$.*

**Corollary 4.10** *Let $(X, d)$ be a bicomplete $T_0$-quasi-pseudometric space, $\phi : X \to \mathbb{R}$ be a bounded from above function and $"\preceq"$ the preorder induced by $\phi$. Let $F : X \times X \to X$ and $G, H : X \to X$ be three $d^s$-sequentially continuous mapping on $X$ such that the pairs $\{F, G\}$ and $\{F, H\}$ satisfy conditions $(C1')$ and $(C2')$. Then, $F, G$ and $H$ have a common coupled fixed point in $X$.*

## Concluding remarks and open problem

All the results given remain true when we replace accordingly the bicomplete quasi-pseudometric space $(X, d)$ by a left Smyth sequentially complete/left $K$-complete or a right Smyth sequentially complete/right $K$-complete space.

Moreover, the reader can convince himself that the proofs are quite straight forward; it is enough to get the right sequence. The major challenge comes when we have more than three maps. Indeed, it is not obvious to see how to construct an appropriate sequence, following the same patent as developed above. More precisely.

Let $(X, d)$ be a Hausdorff left $K$-complete $T_0$-quasi-pseudometric space, $\phi : X \to \mathbb{R}$ be a bounded from above function and $"\preceq"$ the preorder induced by $\phi$. Let $F : X \times X \to X$ and $G_i : X \to X$, $i = 1, 2, \ldots, K$ for $K > 2$ be $K + 1$ $d$-sequentially continuous mapping on $X$ such that the pairs $\{F, G_i\}$, $i = 1, 2, \ldots, K$ are weakly left-related.

(1) Can we prove that $F, G_1, \ldots, G_K$ have a common coupled fixed point in $X$?
(2) Alternatively, what could be a correct formulation of the statement, using the induced preorder and the weakly left-related property that guarantees a positive answer?

## References

1. Agarwal, R.P., El-Gebeily, M.A., O'Regan, D.: Generalized contractions in partially ordered metric spaces. Appl. Anal. **87**(1), 1–8 (2008)
2. Lakshmikantham, V., Ciric, LjB: Coupled fixed point theorems for nonlinear contractions in partially ordered metric space. Nonlinear Anal. **70**(12), 4341–4349 (2009)
3. Mustafa, Z., Sims, B.: Fixed point theorems for contractive mapping in complete G-metric spaces. Fixed Point Theory Appl. (2009). doi:10.1155/2009/917175
4. Bhaskar, T.G., Lakshmikantham, V.: Fixed point theorems in partially ordered metric spaces and applications. Nonlinear Anal. **65**(7), 1379–1393 (2006)
5. Włodarczyk, K., Plebaniak, R.: Asymmetric structures, discontinuous contractions and iterative approximation of fixed and periodic points. Fixed Point Theory Appl. **128**, 1–18 (2013). doi:10.1186/1687-1812-2013-128
6. Gaba, Y.U.: Startpoints and $(\alpha, \gamma)$-contractions in quasi-pseudometric spaces. J. Math. (2014). doi:10.1155/2014/709253 (Article ID 709253, 8 pages)